\documentclass[12pt,twoside]{article}
\usepackage{amsfonts}
{\markboth{\sl J. El Kamel and Ch. Yacoub}
{\sl Almost Everywhere Convergence}\thispagestyle{plain}

\newtheorem{theorem}{Theorem}[section]
\newtheorem{corollary}[theorem]{Corollary}
\newtheorem{lemma}[theorem]{Lemma}
\newtheorem{proposition}[theorem]{Proposition}

\everymath{\displaystyle}

\begin{document}
\title{ Almost Everywhere Convergence of Inverse
Dunkl Transform on the Real Line}
\author{ J. El Kamel and Ch. Yacoub
\\\small Department of Mathematics, Faculty of Sciences of
Monastir,\\ \small 5019 Monastir, TUNISIA
 \\ \small jamel.elkamel@fsm.rnu.tn; chokri.yacoub@fsm.rnu.tn}
\date{}
\maketitle
\begin{abstract}
 In this paper, we will first show that  the maximal operator $S_*^\alpha$ 
 of spherical partial sums $S_R^\alpha$, associated to Dunkl transform on $\mathbb{R}$ is bounded on 
 $L^p\left(\mathbb{R}, \left|x\right|^{2\alpha+1} dx\right)$ functions when 
 $\frac{4(\alpha+1)}{2\alpha+3}<p<\frac{4(\alpha+1)}{2\alpha+1}$, and it implies that, for every 
 $L^p\left(\mathbb{R}, \left|x\right|^{2\alpha+1} dx\right)$ function $f(x)$,  $S_R^\alpha f(x)$
 converges to $f(x)$ almost everywhere as $R\rightarrow \infty$. On the other hand we obtain a 
 sharp version  by showing that $S_*^\alpha$ is bounded from the Lorentz space 
  $L^{p_i,1}\left(\mathbb{R}, \left|x\right|^{2\alpha+1}\right)$ into 
 $L^{p_i,\infty}\left(\mathbb{R}, \left|x\right|^{2\alpha+1}\right),\quad i=0,1$ where 
 $p_0=\frac{4(\alpha+1)}{2\alpha+3}$ and $p_1=\frac{4(\alpha+1)}{2\alpha+1}$.
   
\end{abstract}
{\it Keywords: Dunkl transform, maximal function, almost everywhere convergence, Lorentz space.\\
 
 \begin{section}{Introduction and preliminaries}
 Given $\alpha\geq \frac{-1}{2}$ and a suitable function $f$ on $\mathbb{R}$, its 
 Dunkl transform $D_\alpha$ is defined by
 \begin{equation}
 D_\alpha f(y) = \int_\mathbb{R} f(x) E_\alpha (-ixy)d\mu_\alpha(x), \quad y\in \mathbb{R};
 \end{equation}
 here 
  \begin{equation}
 d\mu_\alpha(x)=\frac{1}{2^{\alpha+1} \Gamma (\alpha+1)}\left|x\right|^{2\alpha+1}dx
 \end{equation}
 \begin{equation}
  E_\alpha (z)=2^\alpha \Gamma (\alpha+1)\left\{ \frac{J_\alpha(iz)}{(iz)^\alpha}+
  z \frac{J_{\alpha+1}(iz)}{(iz)^{\alpha+1}}\right\},
 \end{equation}
 where $J_\alpha$ denotes the Bessel function of the first kind of order $\alpha$. 
 The inverse Dunkl transform $\check{D_\alpha}$ is given by $\check{D_\alpha} f(\lambda)=D_\alpha f(-\lambda).$\\
 \vskip 0.3cm
 In this paper, we are interested in the almost everywhere convergence as $R\rightarrow \infty$ of the 
  partial sums $S_R^\alpha f(x)$ where 
 
 $$S_R^\alpha f(x)= \frac{1}{2^{\alpha+1} \Gamma (\alpha+1)} \int_{\left|y\right|\leq R}
                 D_\alpha f(y) E_\alpha (ixy) \left|y\right|^{2\alpha+1}dy.$$

  Recall that given $\beta\geq -\frac{1}{2}$, the Hankel transform of order $\beta$ of 
   a suitable function $g$ on $(0,\infty)$ is defined by :
   \begin{equation}
   \mathcal{H}_\beta g(y)=\int_0^\infty g(x) \frac{J_\beta(yx)}{(yx)^\beta} x^{2\beta+1}dx,
   \quad y>0.
 \end{equation}
 Nowak and Stempak (\cite{NS}), found an expression of the Dunkl transform $D_\alpha$ in 
  terms of Hankel transform of orders $\alpha$ and $\alpha+1$.\\
 \begin{lemma}(see (\cite{NS})
 Given $\alpha\geq -\frac{1}{2}$, we have :
 \begin{equation}
 D_\alpha f(y)= \mathcal{H}_\alpha (f_e)(\left|y\right|)-
  iy \mathcal{H}_{\alpha+1}\left(\frac{f_o(x)}{x}\right) \left(\left|y\right|\right),
 \end{equation}
 where for a function $f$ on $\mathbb{R}$, we denote by $f_e$ and $f_o$ the restrictions to 
 $(0, \infty)$ of its even and odd parts, respectively, i.e. the functions on $(0, \infty)$ 
  defined by 
  $$ f_e(x)=\frac{1}{2}\left(f(x)+f(-x)\right),\quad f_o(x)=\frac{1}{2}\left(f(x)-f(-x)\right),\quad x>0.$$
 \end{lemma}
 Define, the partial sums $s_R^\beta g(x)$ by :
 \begin{equation}
 s_R^\beta g(x)=\int_0^R \mathcal{H}_\beta g(y) \frac{J_\beta(xy)}{(xy)^\beta} y^{2\beta+1}dy,
 \quad x>0
 \end{equation}
 and 
 \begin{equation}
 s_*^\beta g(x)=\sup_{R>0}\left|s_R^\beta g(x)\right|.
 \end{equation}
 In 1988, Y. Kanjin (\cite{K}) and E. Prestini (\cite{P}) proved, independently, the following :
 \begin{theorem}
 Let $\beta\geq -\frac{1}{2}$.\\
 $\bullet$ If \quad $\frac{4(\beta+1)}{2\beta+3}<p<\frac{4(\beta+1)}{2\beta+1}$ then $s_*^\beta $ is bounded 
 on $L^p\left((0,\infty),x^{2\beta+1}\right)$ functions.\\
 $\bullet$ If \quad $p\leq \frac{4(\beta+1)}{2\beta+3}$ or  $p\geq \frac{4(\beta+1)}{2\beta+1}$ then
  $s_*^\beta $ is not bounded on $L^p\left((0,\infty),x^{2\beta+1}\right)$ functions.
 \end{theorem}
 Throughout this paper we use the convention that $c_{\alpha}$ denotes a constant, depending on 
 $\alpha$ and $p$, its value may change from line to line.
 \end{section}
 \begin{section}{Almost everywhere convergence}
 
 Define linear operators $S_R^\alpha, R>0$ and $S_*^\alpha$ on the Schwartz space 
 $S\left(\mathbb{R}\right)$ by
 \begin{equation} 
 S_R^\alpha f(x)= \frac{1}{2^{\alpha+1} \Gamma (\alpha+1)} \int_{\left|y\right|\leq R}
                 D_\alpha f(y) E_\alpha (ixy) \left|y\right|^{2\alpha+1}dy
\end{equation} 
and              
 \begin{equation} 
 S_*^\alpha f(x)=\sup_{R>0}\left|S_R^\alpha f(x)\right|, \quad x\in \mathbb{R}.
 \end{equation}
 \begin{lemma}
 Given $\alpha\geq -\frac{1}{2}$, we have
 \begin{equation} 
 S_R^\alpha(f)(x)=s_R^\alpha(f_e)(\left|x\right|)+ x s_R^{\alpha+1}\left(\frac{f_o(r)}{r}\right)
 (\left|x\right|),
 \end{equation} 
 \begin{equation} 
 S_*^\alpha f(x)\leq s_*^\alpha(f_e)(\left|x\right|)+ \left|x\right|s_*^{\alpha+1}
 \left(\frac{f_o(r)}{r}\right)(\left|x\right|).
 \end{equation} 
  \end{lemma}
 \begin{proof}
 Let $x\in \mathbb{R}$. By (3), (8) and lemma 1.1, we have\\
 
  $S_R^\alpha f(x) =\frac{1}{2^{\alpha+1} \Gamma (\alpha+1)} \int_{\left|y\right|\leq R}\left[\mathcal{H}_\alpha (f_e)(\left|y\right|)-
  iy \mathcal{H}_{\alpha+1}
  \left(\frac{f_o(r)}{r}\right) \left(\left|y\right|\right)\right]$
 $$\qquad\qquad\qquad\qquad \qquad\qquad\qquad\left[2^\alpha \Gamma(\alpha+1)
 \left\{\frac{J_\alpha (yx)}{(yx)^\alpha}+
 ixy \frac{J_{\alpha+1} (yx)}{(yx)^{\alpha+1}}\right\}\right]\left|y\right|^{2\alpha+1}dy
 $$
 $$=\frac{1}{2} \int_{\left|y\right|\leq R}\mathcal{H}_\alpha (f_e)(\left|y\right|)
 \frac{J_\alpha (yx)}{(yx)^\alpha}\left|y\right|^{2\alpha+1}dy\qquad\qquad\qquad$$
 $$+\frac{ix}{2}
 \int_{\left|y\right|\leq R}y\mathcal{H}_\alpha (f_e)(\left|y\right|)
 \frac{J_{\alpha+1} (yx)}{(yx)^{\alpha+1}}\left|y\right|^{2\alpha+1}dy\qquad\qquad$$
 $$\qquad\qquad\quad -\frac{i}{2}\int_{\left|y\right|\leq R}
 y \mathcal{H}_{\alpha+1}\left(\frac{f_o(r}{r}\right) \left(\left|y\right|\right)
 \frac{J_\alpha (yx)}{(yx)^\alpha}
 \left|y\right|^{2\alpha+1}dy \qquad\qquad\qquad\qquad$$
 $$\qquad\qquad\quad +\frac{x}{2} 
 \int_{\left|y\right|\leq R} \mathcal{H}_{\alpha+1}\left(\frac{f_o(r)}{r}\right) 
 \left(\left|y\right|\right)\frac{J_{\alpha+1} (yx)}{(yx)^{\alpha+1}}
 \left|y\right|^{2\alpha+3}dy\qquad\qquad\qquad\qquad$$
 \noindent We note that  the second and the third integrals are equal to zero. So
 
 \noindent$S_R^\alpha f(x)=\frac{1}{2} \int_{\left|y\right|\leq R}\mathcal{H}_\alpha (f_e)(\left|y\right|)
 \frac{J_\alpha (yx)}{(yx)^\alpha}\left|y\right|^{2\alpha+1}dy$
 $$+
 \frac{x}{2} 
 \int_{\left|y\right|\leq R} \mathcal{H}_{\alpha+1}\left(\frac{f_o(r)}{r}\right) 
 \left(\left|y\right|\right)\frac{J_{\alpha+1} (yx)}{(yx)^{\alpha+1}}.
 \left|y\right|^{2\alpha+3}dy\qquad\qquad$$
 
 $$\qquad\qquad=\int_0^R \mathcal{H}_\alpha (f_e)(y)
 \frac{J_\alpha (\left|x\right|y)}{(\left|x\right|y)^\alpha} y^{2\alpha+1}dy +
 x\int_0^R \mathcal{H}_{\alpha+1}\left(\frac{f_o(r)}{r}\right) 
 \left(y\right)\frac{J_{\alpha+1} (\left|x\right|y)}{(\left|x\right|y)^{\alpha+1}}
 y^{2\alpha+3}dy$$
 $$=s_R^\alpha(f_e)(\left|x\right|)+ x s_R^{\alpha+1}\left(\frac{f_o(r)}{r}\right)
 (\left|x\right|).\qquad\qquad\qquad\qquad\qquad$$
 Thus
 $$S_*^\alpha f(x)\leq s_*^\alpha(f_e)(\left|x\right|)+ \left|x\right|s_*^{\alpha+1}
 \left(\frac{f_o(r)}{r}\right)(\left|x\right|).$$ 
 \end{proof}
 
 \begin{proposition}
 Let $\alpha\geq -\frac{1}{2}$.\\
  $\bullet$ If \quad $\frac{4(\alpha+1)}{2\alpha+3}<p<\frac{4(\alpha+1)}{2\alpha+1}$ then $S_*^\alpha $ is bounded 
 on $L^p\left(\mathbb{R},\left|x\right|^{2\alpha+1}dx\right)$ functions.\\
 $\bullet$ If \quad $p\leq \frac{4(\alpha+1)}{2\alpha+3}$ or  $p\geq \frac{4(\alpha+1)}{2\alpha+1}$ then
  $S_*^\alpha$ is not bounded on $L^p\left(\mathbb{R},\left|x\right|^{2\alpha+1}dx\right)$ functions.
 \end{proposition}
 \begin{proof}
 $S_*^\alpha $ cannot be bounded for $p\leq\frac{4(\alpha+1)}{2\alpha+3}$ or 
 $p\geq \frac{4(\alpha+1)}{2\alpha+1}$ (see: \cite{K}, \cite{P}).\\
 By theorem 1, we have for $\frac{4(\alpha+1)}{2\alpha+3}<p<\frac{4(\alpha+1)}{2\alpha+1}$ 
 $$\left\|s_*^\alpha(f_e)(\left|x\right|)\right\|_{L^p\left(\mathbb{R},\left|x\right|^{2\alpha+1}dx\right)}
 =2 \left\|s_*^\alpha(f_e)\right\|_{L^p\left((0,\infty),x^{2\alpha+1}dx\right)}\qquad \qquad\qquad\qquad\qquad$$
  $$\leq  c_\alpha\left\|f_e\right\|_{L^p\left((0,\infty),x^{2\alpha+1}dx\right)}$$
  $$\leq c_\alpha \left\|f\right\|_{L^p\left(\mathbb{R},\left|x\right|^{2\alpha+1}dx\right)}.\quad$$
 On the other hand, as in (\cite{P},\cite{RS}), one gets
 \begin{equation}
 \left|x\right|s_*^{\alpha+1}
 \left(\frac{f_o(r)}{r}\right)(\left|x\right|)\leq 
 \frac{c_\alpha}{\left|x\right|^{\alpha+\frac{1}{2}}} \left[M+H+\widetilde{H}+\widetilde{C}\right]\left[\frac{f_o(r)}{r} 
 r^{\alpha+\frac{3}{2}}\right](\left|x\right|),
 \end{equation}
 where $M,H,\widetilde{H}$ and $\widetilde{C}$ denotes respectively, the maximal 
 function, the Hilbert integral, the maximal Hilbert transform and the Carleson operator.\\
 Let $K=M+H+\widetilde{H}+\widetilde{C}$ and $w\in A_p\left(\mathbb{R}\right), p>1$. It is well known that
 \begin{equation}
 \left\|Kf\right\|_{L^p\left(\mathbb{R}, w(x) dx\right)}\leq c_\alpha 
 \left\|f\right\|_{L^p\left(\mathbb{R}, w(x) dx\right)}.
 \end{equation}
 Hence\\
 
  $\left\|\left|x\right|s_*^{\alpha+1}
 \left(\frac{f_o(r)}{r}\right)(\left|x\right|)\right\|_
 {L^p\left(\mathbb{R}, \left|x\right|^{2\alpha+1} dx\right)}
 $
 $$\qquad\qquad\qquad\qquad\qquad\leq
 c_{\alpha}\left\|\left|x\right|^{-\alpha-\frac{1}{2}} K
 \left[\frac{f_o(r)}{r}r^{\alpha+\frac{3}{2}}\right](\left|x\right|)\right\|_
 {L^p\left(\mathbb{R}, \left|x\right|^{2\alpha+1} dx\right)}$$
 $$\qquad\qquad\qquad\leq c_\alpha \left\|K\left[\frac{f_o(r)}{r}r^{\alpha+3/2}\right](\left|x\right|)\right\|_
 {L^p\left(\mathbb{R}, w(x) dx\right)},$$
 with $w(x)=\left|x\right|^{2\alpha+1-p(\alpha+1/2)}$.\\

\noindent Since $\frac{4(\alpha+1)}{2\alpha+3}<p<\frac{4(\alpha+1)}{2\alpha+1}$
  if and only if $-1< 2\alpha+1-p(\alpha+1/2)<p-1$,\\ 
  
\noindent  then $w\in A_p\left(\mathbb{R}\right)$ and by (13)

 $$\left\|\left|x\right|s_*^{\alpha+1}
 \left(\frac{f_o(r)}{r}\right)(\left|x\right|)\right\|_
 {L^p\left(\mathbb{R}, \left|x\right|^{2\alpha+1} dx\right)}\leq c_\alpha
 \left\|\frac{f_0(\left|x\right|)}{\left|x\right|}\left|x\right|^{\alpha+3/2}\right\|_
 {L^p\left(\mathbb{R}, w(x) dx\right)}$$
 $$\qquad\qquad\qquad\qquad\qquad\qquad
 \leq c_{\alpha}\left\|f_o(x)\right\|_{L^p\left(\mathbb{R}, \left|x\right|^{2\alpha+1} dx\right)}$$
 $$\qquad\qquad\qquad\qquad\qquad\qquad
 \leq c_{\alpha}\left\|f(x)\right\|_{L^p\left(\mathbb{R}, \left|x\right|^{2\alpha+1} dx\right)}.$$
 We conclude by lemma 2.1.
 \end{proof}
 \begin{corollary}
 For every $f\in L^p\left(\mathbb{R}, \left|x\right|^{2\alpha+1} dx\right)$,
 if\quad
 $\frac{4(\alpha+1)}{2\alpha+3}<p<\frac{4(\alpha+1)}{2\alpha+1}$,\quad then
 $$S_R^\alpha f(x)\rightarrow f(x) \quad a.e.\quad as \quad R\rightarrow \infty$$
 \end{corollary}
\end{section}
\begin{section}{Endpoint estimates} 
 We recall that the Lorentz space $L^{p,q}\left(X,\mu\right)$, , 
 is the set of all measurable functions $f$ on $X$ satisfying 
 $$\left\|f\right\|_{p,q}=\left(\frac{q}{p}\int_0^\infty \left(t^{\frac{1}{p}}f^*(t)\right)^q \frac{dt}{t}\right)^{\frac{1}{q}} <\infty$$
 when
 $1\leq p<\infty$, $1\leq q< \infty$, and 
 $$\left\|f\right\|_{p,q}=\sup_{t>0} t^{\frac{1}{p}}f^*(t)=
 \sup_{\lambda>0}\lambda \left(d_f(\lambda)\right)^{\frac{1}{p}}<\infty$$
 when
 $1\leq p\leq \infty$ and  $q= \infty$. Where $f^*$ denotes the nonincreasing rearrangement 
 of $f$, i.e. 
 $$f^*(t)= \inf\left\{s>0/ d_f(s)\leq t\right\}, 
 \qquad d_f(s)=\mu\left\{x\in X/ \left|f(x)\right|>s\right\}.$$
 In  1991 E. Romera and F. Soria \cite{RS} (see also L. Colzani and all \cite{CCTV}) proved 
 the following :
  
 \begin{theorem}
 Let $\alpha > -\frac{1}{2}$, then $s_*^\alpha$ is bounded from the Lorentz space 
 $L^{p_i,1}\left((0,\infty), x^{2\alpha+1}dx\right)$ into 
 $L^{p_i,\infty}\left((0,\infty), x^{2\alpha+1}dx\right)$, i=0,1 when 
 $p_0= \frac{4(\alpha+1)}{2\alpha+3}$ and $p_1= \frac{4(\alpha+1)}{2\alpha+1}$ is 
 the index conjugate to $p_0$.
 \end{theorem}
Using this result, we will see that proposition 2.2 can be strengthened. More precisely 
 we obtain :
 \begin{proposition}
 
 Let $\alpha > -\frac{1}{2}$, then $S_*^\alpha$ is bounded from the Lorentz space 
 $L^{p_i,1}\left(\mathbb{R}, \left|x\right|^{2\alpha+1}dx\right)$ into 
 $L^{p_i,\infty}\left(\mathbb{R}, \left|x\right|^{2\alpha+1}dx\right), i=0,1.$
 
 \end{proposition}
 So using the formulation of Marcinkiewicz interpolation theorem in terms of 
  Lorentz space we retrieve Proposition 2.2 $(\alpha>-\frac{1}{2})$ as a corollary.\\
 \begin{proof}
 By lemma 2.1, we have \\
 
 $\mu_\alpha\left\{x\in\mathbb{R} / S_*^\alpha f(x) > \lambda\right\}\leq
  \mu_\alpha\left\{x\in\mathbb{R} / s_*^\alpha f_e(\left|x\right|) > \frac{\lambda}{2}\right\}$
  $$\qquad\qquad\qquad\qquad \qquad\quad+ \mu_\alpha\left\{x\in\mathbb{R} / \left|x\right|
 s_*^{\alpha+1} \left(\frac{f_o(r)}{r}\right)(\left|x\right|) > \frac{\lambda}{2}\right\}.$$
 $$= I+II\quad\qquad\qquad$$
 By theorem 2.4, we get :\\
 
 $\mu_\alpha\left\{x\in\mathbb{R} / s_*^\alpha f_e(\left|x\right|) > \frac{\lambda}{2}\right\}
 =2\mu_\alpha\left\{x\in (0,\infty) / s_*^\alpha f_e(x) > \frac{\lambda}{2}\right\}$
 $$\qquad \qquad\qquad\leq \frac{c_\alpha}{\lambda^{p_i}} \left\|f_e\right\|_{p_i,1}
 \leq \frac{c_\alpha}{\lambda^{p_i}} \left\|f\right\|_{p_i,1}.$$
 
 To estimate $II$, we follow closely \cite{RS} and we sketch a proof for completeness. We 
  decompose the set :
  
 $\displaystyle\left\{x \in \mathbb{R}/ 
 \left|x\right|s_*^{\alpha+1} \left(\frac{f_o(r)}{r}\right)(\left|x\right|) >\frac{\lambda}{2}\right\}$
 $$\qquad\qquad\qquad\qquad=\bigcup_{k\in \mathbb{Z}} \left\{x \in \mathbb{R}/ 
 \left|x\right|\in I_k, 
 \left|x\right|s_*^{\alpha+1} \left(\frac{f_o(r)}{r}\right)(\left|x\right|) > \frac{\lambda}{2}\right\},$$
 where $I_k=[2^k,2^{k+1}[$.\\
 
 \noindent Put $g(r):= \frac{f_o(r)}{r}=g_k^1(r)+g_k^2(r)$, with $g_k^1=g\chi_{I_k^*}$, $g_k^2=g\chi_{(I_k^*)^c}$
 , where $I_k^*=]2^{k-1},2^{k+2}[$.\\
 
 \noindent By (12), we have :
 $$
 \left|x\right|s_*^{\alpha+1} \left(g_k^1(r)\right)(\left|x\right|)\leq 
 \frac{c_\alpha}{\left|x\right|^{\alpha+1/2}} K\left(g_k^1(r) r^{\alpha+3/2}\right)(\left|x\right|).
  $$
 By (\cite{RS}, p: 1021), we have for $1<p<\infty$, \\
 $\sum_{k\in\mathbb{Z}}
 \mu_\alpha\left\{ x\in \mathbb{R} / 
 \left|x\right|\in I_k, \frac{1}{\left|x\right|^{\alpha+1/2}} 
 K\left(g_k^1(r) r^{\alpha+3/2}\right)(\left|x\right|)>\frac{\lambda}{2} \right\}$
 $$\leq \frac{c_{\alpha}}{\lambda^p} 
 \left\|f_o\right\|_{L^p\left(\mathbb{R},\left|x\right|^{2\alpha+1}dx\right)}^p
  \leq \frac{c_{\alpha}}{\lambda^p} 
 \left\|f\right\|_{L^p\left(\mathbb{R},\left|x\right|^{2\alpha+1}dx\right)}^p
  \leq \frac{c_{\alpha} }{\lambda^p} \left\|f\right\|_{p,1}^p.$$
 
\noindent On the other hand as in (\cite{RS}, p: 1021), we have \\

 $\left|x\right|s_*^{\alpha+1} \left(g_k^2(r)\right)(\left|x\right|)\leq 
 \frac{c_\alpha}{\left|x\right|^{\alpha+1/2}} \int_0^\infty 
 \frac{s^{\alpha+3/2}\left|f_o(s)\right|}{s(\left|x\right|+s)} ds$
 $$\quad\leq 
 \frac{c_\alpha}{\left|x\right|^{\alpha+3/2}} \int_0^\infty \left|f_o(s)\right| s^{\alpha+1/2} ds$$
 $$\qquad\qquad\leq 
 \frac{c_\alpha}{\left|x\right|^{\alpha+3/2}} \int_\mathbb{R} \left|f_o(s)\right| 
 \frac{1}{\left|s\right|^{\alpha+1/2}} \left|s\right|^{2\alpha+1} ds.$$
 Remark that we have considered $f_0$ as a function defined on $\mathbb{R}$.\\
\noindent As the same we get \\
 $\left|x\right|s_*^{\alpha+1} \left(g_k^2(r)\right)(\left|x\right|)\leq 
 \frac{c_\alpha}{\left|x\right|^{\alpha+1/2}} \int_0^\infty 
 \left|f_o(s)\right| s^{\alpha-1/2} ds$
 $$\qquad \leq\frac{c_\alpha}{2\left|x\right|^{\alpha+1/2}} \int_\mathbb{R} \left|f_o(s)\right| \frac{1}{\left|s\right|^{\alpha+3/2}} \left|s\right|^{2\alpha+1} ds.$$\\
 Using the following facts :
 $$
 \frac{1}{\left|x\right|^{\alpha+\frac{1}{2}}} \in L^{p_1,\infty}\left(\mathbb{R}, \left|x\right|^{2\alpha+1}\right),
 $$
 
 $$ \frac{1}{\left|x\right|^{\alpha+\frac{3}{2}}} \in L^{p_0,\infty}\left(\mathbb{R}, \left|x\right|^{2\alpha+1}\right),
 $$
and Holder's inequality for the Lorentz spaces, we arrive to :
$$
 \mu_\alpha\left\{x\in\mathbb{R} / \left|x\right|s_*^{\alpha+1} 
 \left(g_k^2(r)\right)(\left|x\right|)> \frac{\lambda}{2}\right\}\leq \frac{c_\alpha}{\lambda^{p_i}}
 \left\|f_o\right\|_{p_i,1}^{p_i}\leq 
 \frac{c_\alpha}{\lambda^{p_i}}
 \left\|f\right\|_{p_i,1}^{p_i},
 $$
 which completes the proof.
 
 \end{proof}
\vskip 0.5 cm
\noindent{\bf Acknowledgment.} We are grateful to Professor K. Stempak for sending us the preprint
 \cite{NS}.
\end{section}

 \end{document}